\documentclass[11pt,a4paper]{article}
\usepackage[T1]{fontenc}
\usepackage{lmodern,microtype}
\usepackage{amsmath,amssymb,amsthm,mathtools}
\usepackage{geometry,booktabs,array,longtable,enumitem,xcolor}
\usepackage[hidelinks]{hyperref}
\newtheorem{theorem}{Theorem}[section]
\newtheorem{lemma}[theorem]{Lemma}
\newtheorem{corollary}[theorem]{Corollary}
\newtheorem{proposition}[theorem]{Proposition}

\theoremstyle{definition}\newtheorem{definition}[theorem]{Definition}
\theoremstyle{remark}
\newcommand{\CN}{\mathrm{CN}}
\newcommand{\ON}{\mathrm{ON}}
\newcommand{\rvc}{\operatorname{rvc}}
\newcommand{\srvc}{\operatorname{srvc}}

\newcommand{\diam}{\operatorname{diam}}

\title{Coloring Variants on Split, Starlike and
Well-Partitioned Chordal Graphs:
Algorithms, Bounds and Hardness}

\author{
Suchanda Roy,
Ramesh Hariharasubramanian\thanks{Corresponding author},
Tithi Dwary,
Biswajit Das\\[6pt]
\small Department of Mathematics,
Indian Institute of Technology Guwahati\\
\small Guwahati, Assam 781039, India\\[4pt]
\small 
\texttt{r.suchanda@iitg.ac.in},\\
\texttt{ramesh\_h@iitg.ac.in},\\
\texttt{tithi.dwary@iitg.ac.in},\\
\texttt{d.biswajit@iitg.ac.in}
}

\date{}

\begin{document}
\maketitle

\begin{abstract}
In this paper, we study full and partial conflict-free
coloring of open and closed neighborhoods, together with
rainbow vertex connection and strong rainbow vertex
connection, on circular-convex split graphs, starlike
graphs, and well-partitioned chordal (WPC) graphs.
Using the known polynomial-time algorithm for
conflict-free coloring of interval hypergraphs, we provide
polynomial-time computable upper bounds for
conflict-free open-neighborhood coloring on
circular-convex split and starlike graphs.
We further extend these bounds to circular-convex
WPC graphs by using partition tree.

For closed-neighborhood conflict-free coloring,
we obtain an exact characterization of both the partial
and full variant on starlike graphs and establish upper
bounds of two and three colors for partial and
full coloring, respectively, on WPC graphs.
We also prove NP-completeness of deciding whether
a WPC graph admits a one-color partial or a
two-color full closed-neighborhood conflict-free
coloring.

For rainbow vertex connection, we obtain exact
results on circular-convex split and starlike
graphs and investigate their extensions to WPC
graphs through split components and partition
trees. We establish component-based upper bounds,
improve them by reusing colors across certain
leaf bags, and identify subclasses of WPC graphs
for which these bounds are attained.
\end{abstract}


\section{Introduction}
\label{sec:introduction}

Graph coloring is a fundamental problem in graph theory
with applications in communication networks, resource
allocation, scheduling, and information transmission.
Several variants of graph coloring have been introduced
to study different structural and optimization
requirements. Among these, conflict-free coloring and
rainbow vertex coloring have received considerable
attention from both structural and algorithmic
perspectives.

In a conflict-free coloring, the objective is to assign
colors to the vertices of a graph such that every vertex
has a color appearing exactly once in its neighborhood.
Depending on whether open or closed neighborhoods are
considered, we obtain conflict-free open-neighborhood
coloring (CFON) and conflict-free closed-neighborhood
coloring (CFCN), respectively. In the full coloring
variants, every vertex must receive a color, whereas
in the partial variants, some vertices may remain
uncolored. These four parameters behave differently,
and their computational complexities may vary even
on the same graph class.

Conflict-free coloring has been studied on several
graph classes, including interval graphs, split graphs,
bipartite graphs, and chordal graphs. Reddy
\cite{Reddy2018} investigated conflict-free coloring
problems on restricted graph classes and established
several complexity results. More recently, Bhyravarapu
et al.~\cite{BKM2025} studied conflict-free coloring
on subclasses of perfect and bipartite graphs.
On the other hand, Dhannya and Narayanaswamy
\cite{Dhannya2020} developed a polynomial-time
algorithm for conflict-free coloring of interval
hypergraphs. This result provides an important
algorithmic tool for studying conflict-free coloring
on graph classes whose neighborhoods satisfy
interval or convexity properties.

In this paper, we investigate how such interval
structures can be used to construct conflict-free
colorings on circular-convex split graphs and
starlike graphs. A split graph is a graph whose
vertex set can be partitioned into a clique $C$
and an independent set $I$. It is circular-convex
with respect to $C$ if the vertices of $C$ admit
a circular ordering in which the neighborhood of
every vertex of $I$ forms a circular interval.
This property allows us to represent the
neighborhoods of independent vertices as
hyperedges on a circle. However, satisfying
these hyperedges alone does not necessarily
satisfy the conflict-free coloring requirements
of the clique vertices. Therefore, additional
conditions are needed to obtain a valid coloring
of the entire graph.

We also consider well-partitioned chordal (WPC)
graphs, introduced and studied by Ahn et
al.~\cite{Ahn2022}. These graphs admit a partition
of their vertices into cliques arranged along
a tree, where the edges between adjacent cliques
form complete bipartite graphs between specified
boundary sets. This partition-tree structure
provides a natural framework for extending
coloring techniques from split graphs to more
general graph classes. In particular, we study
bagwise circular-convex WPC graphs, where the
boundary sets of every bag form circular intervals
in a suitable ordering.

Another important coloring problem considered
in this paper is rainbow vertex connection.
A vertex-colored graph is rainbow vertex-connected
if every pair of vertices is joined by a path
whose internal vertices have pairwise distinct
colors. If such a path is additionally required
to be a shortest path, the coloring is called
strong rainbow vertex-connected. The corresponding
minimum numbers of colors are denoted by
$\rvc(G)$ and $\srvc(G)$, respectively.
These problems have been investigated on several
restricted graph classes. Heggernes et
al.~\cite{Heggernes2018} established hardness
results on split graphs, while Lima, van Leeuwen,
and van der Wegen~\cite{Lima2021} developed
polynomial-time algorithms for rainbow vertex
coloring on certain subclasses of chordal graphs.

Motivated by these results, we study conflict-free
coloring and rainbow vertex connection on
circular-convex split, starlike, and WPC graphs.
Our main contributions are summarized below.

\begin{itemize}

\item We develop anchored interval-hypergraph
coloring constructions that provide polynomial-time
computable upper bounds for CFON and its partial
variant on circular-convex split and starlike
graphs. We also construct graph families showing
that the logarithmic order of these bounds is tight.

\item We extend the CFON coloring constructions
to bagwise circular-convex WPC graphs using three
color palettes based on the depths of the bags
in a partition tree. This gives polynomial-time
constructive upper bounds for both full and
partial CFON coloring.

\item For CFCN, we obtain an exact characterization
for both the full and partial conflict-free chromatic number
on starlike graphs. We also establish upper
bounds of two and three colors for partial
and full CFCN, respectively, on WPC graphs.
Furthermore, we prove NP-completeness of deciding
whether a WPC graph admits a one-color partial
CFCN coloring or a two-color full CFCN coloring. 

\item For rainbow vertex connection, we study
the relationship between $\rvc(G)$ and $\srvc(G)$
on split graphs and obtain exact expressions
for circular-convex split and starlike graphs.
We further investigate these parameters on
WPC graphs using their split components and
partition-tree structures.
We establish component-based upper bounds
for rainbow vertex connection on WPC graphs
and improve them by reusing colors among certain
leaf bags. We also identify subclasses of WPC
graphs for which these bounds are exact.

\end{itemize}

The remainder of the paper is organized as
follows. Section~\ref{sec:preliminaries}
introduces the necessary definitions,
notations, and previously established results.
Section~\ref{sec:cfon} presents our results
on CFON and its partial variant for convex
split, starlike, and WPC graphs.
Section~\ref{sec:cfcn} studies CFCN and
its partial variant, including the
NP-completeness results on WPC graphs.
Section~\ref{sec:rainbow} investigates
rainbow vertex connection and strong
rainbow vertex connection on the
considered graph classes.
Finally, Section~\ref{sec:conclusion}
concludes the paper and discusses
directions for future research.

\section{Preliminaries}
\label{sec:preliminaries}

All graphs considered in this paper are finite,
simple, and undirected. For a graph $G=(V,E)$
and a vertex $v\in V(G)$, let $N(v)$ and
$N[v]=N(v)\cup\{v\}$ denote the open and closed
neighborhoods of $v$, respectively.
We use $\diam(G)$ to denote the diameter of
a connected graph $G$. Throughout the paper,
graphs considered for conflict-free
open-neighborhood coloring have no isolated
vertices.

A graph $G$ is called a \emph{split graph} if
its vertex set can be partitioned into a clique
$C$ and an independent set $I$. A split graph
is \emph{circular-convex with respect to $C$}
if there exists a cyclic ordering of the
vertices of $C$ such that $N(x)\cap C$
forms a circular interval for every
$x\in I$. In particular, a circular interval
consists of consecutive vertices in the
cyclic ordering.

A \emph{starlike graph} consists of a central
clique $C$ and pairwise anticomplete peripheral
cliques $Q_1,\ldots,Q_m$, where all vertices
of each $Q_i$ have the same nonempty
neighborhood $S_i\subseteq C$. We call $S_i$
the \emph{boundary} of $Q_i$.
When all these boundaries are circular
intervals in a common cyclic ordering of $C$,
we call the graph circular-convex starlike.

We first recall the definitions of the
conflict-free coloring parameters studied
in this paper.

\begin{definition}
A \emph{full conflict-free open-neighborhood
(CFON) $k$-coloring} of a graph $G$ is a
mapping $c:V(G)\to[k]$ such that, for every
$v\in V(G)$, some color occurs exactly once
in $N(v)$. The minimum such $k$ is denoted
by $\chi_{\ON}(G)$.

In a \emph{partial CFON $k$-coloring},
the mapping is
$c:V(G)\to\{0,1,\ldots,k\}$, where color
$0$ represents an uncolored vertex and
cannot serve as a witnessing color.
The corresponding minimum number of
positive colors is denoted by
$\chi_{\ON}^*(G)$.

Replacing the open neighborhood $N(v)$
with the closed neighborhood $N[v]$
defines \emph{full conflict-free
closed-neighborhood coloring (CFCN)}
and \emph{partial CFCN}, whose minimum
numbers of colors are denoted by
$\chi_{\CN}(G)$ and
$\chi_{\CN}^*(G)$, respectively.
\end{definition}

Since every full coloring is also a
partial coloring, and uncolored vertices
can be assigned one additional color,
we have
$
\chi_{\ON}^*(G)\leq\chi_{\ON}(G)
\leq\chi_{\ON}^*(G)+1
$
and
$
\chi_{\CN}^*(G)\leq\chi_{\CN}(G)
\leq\chi_{\CN}^*(G)+1.
$
Thus, although the full and partial
parameters differ by at most one,
their exact determination may require
different arguments.

Next, we recall the structure of
well-partitioned chordal graphs,
which plays an important role in
extending our coloring constructions.

\begin{definition}\cite{Ahn2022}
A connected graph $G$ is
\emph{well-partitioned chordal (WPC)}
if its vertex set admits a partition
into nonempty clique bags $\mathcal B$,
indexed by the vertices of a tree $T$,
such that no edges exist between
nonadjacent bags. Furthermore, for
every pair of adjacent bags $B,D$,
the edges between them form a
complete bipartite graph between
nonempty boundary sets
$S_{B,D}\subseteq B$ and
$S_{D,B}\subseteq D$.

For each bag $B$, we associate a
\emph{split component} $H_B$ whose
clique is $B$ and whose independent
set contains one marker vertex
$m_{B,D}$ for every
$D\in N_T(B)$, with
$
N_{H_B}(m_{B,D})=S_{B,D}.
$
The representation is called
\emph{bagwise circular-convex}
if all boundary sets of each bag
$B$ are circular intervals in a
common cyclic ordering of $B$.
\end{definition}

We now recall a known algorithmic
result for conflict-free coloring
of interval hypergraphs. This result
will be used as a subroutine in
our coloring constructions.

\begin{theorem}\cite{Dhannya2020}
Given an interval hypergraph
together with its linear interval
representation, an optimal partial
conflict-free coloring of its
hyperedges can be computed in
polynomial time. Moreover, full
interval-hypergraph coloring can
be reduced to the partial variant
by adding all singleton hyperedges.
\end{theorem}

We also use the connection between
partial closed-neighborhood
conflict-free coloring and efficient
domination. The notion of efficient
domination was studied by Bange,
Barkauskas, and Slater
\cite{BangeBarkauskasSlater},
and its connection with conflict-free
coloring is discussed in
\cite{AbelEtAl2018}.

\begin{lemma}\cite{AbelEtAl2018}]
\label{lem:effdom}
For every nonempty graph $G$,
$\chi_{\CN}^*(G)=1$ if and only if
$G$ has an efficient dominating
set $D$, that is,
$
|N[v]\cap D|=1
\qquad\text{for every }v\in V(G).
$
\end{lemma}

Indeed, assigning color $1$ exactly
to the vertices of an efficient
dominating set gives a one-color
partial CFCN coloring, and the
converse follows directly from
the definition.

Finally, we recall the notions of
rainbow vertex connection and
strong rainbow vertex connection.

A vertex-colored path is called a
\emph{rainbow vertex path} if all
its internal vertices have pairwise
distinct colors. A graph $G$ is
\emph{rainbow vertex-connected}
if every pair of vertices is joined
by a rainbow vertex path. The
minimum number of colors required
is denoted by $\rvc(G)$. If every pair of vertices is joined
by a rainbow vertex path that is
also a shortest path, then $G$
is called \emph{strong rainbow
vertex-connected}. The corresponding
minimum number of colors is denoted
by $\srvc(G)$. 

For complete graphs, we adopt the
convention
$\rvc(G)=\srvc(G)=0$.
For every connected graph $G$,
the following inequalities hold:
$
\diam(G)-1
\leq\rvc(G)
\leq\srvc(G).
$ 

Rainbow vertex coloring is
computationally hard even on split
graphs~\cite{Heggernes2018}. On the other hand, Lima,
van Leeuwen, and van der
Wegen~\cite{Lima2021} developed
polynomial-time algorithms for
rainbow vertex connection and
strong rainbow vertex connection
on split strongly chordal graphs.

These definitions and known results
form the basis for the coloring
bounds, structural characterizations,
and complexity results presented
in the following sections.

\section{Interval hypergraphs and CFON and CFON* on Convex Classes}\label{sec:cfon}
Let $(X,\mathcal S)$ be a family of nonempty circular intervals on a cyclically ordered set $X$, $|X|\ge2$. If $H=(U,\mathcal F)$ is a hypergraph, let $\mathrm{scf}_1(H)$ be the minimum number of colors required for a partial conflict-free coloring of its hyperedges, and let $\mathrm{cf}(H)$ denote the full optimum. By adding singleton hyperedges, full interval-hypergraph coloring reduces to the partial problem; the interval-hypergraph optimum is polynomial-time computable \cite{Dhannya2020}.

For $p\in X$ we let $X_p=X\setminus\{p\}$ and
$
 H_p(X,\mathcal S)=\left(X_p,\{S\in\mathcal S:p\notin S\}\cup\{X_p\}\right),
$
and define
\[
 a^*(X,\mathcal S)=1+\min_{p\in X}\mathrm{scf}_1(H_p),\qquad
 a(X,\mathcal S)=1+\min_{p\in X}\mathrm{cf}(H_p).
\]
When $|X|=1$, we set both parameters to one. Removing $p$ converts all hyperedges of $H_p$ into linear intervals.

\begin{lemma}\label{lem:anchor}
Both $a^*(X,\mathcal S)$ and $a(X,\mathcal S)$ are computable in polynomial time from a cyclic representation. There exist respective partial and full colorings of $X$ using these numbers of colors such that every $S\in\mathcal S$ contains a uniquely occurring positive color and two distinct colors each occur exactly once on $X$. Furthermore,
$
 a^*(X,\mathcal S)\le a(X,\mathcal S)\le1+\lceil\log_2|X|\rceil.
$
\end{lemma}
\begin{proof}
For a fixed anchor $p$, we optimally color $H_p$ and give $p$ a fresh color. The additional hyperedge $X_p$ ensures a positive color $\beta$ occurs exactly once on $X_p$ at some $q$; the fresh color at $p$ and $\beta$ occur exactly once on $X$. Each interval $S$ either contains $p$ and uses its fresh color, or belongs to $H_p$ and has an existing witness. we enumerate $p$. A recursive midpoint coloring of the linearly ordered $|X|-1$ remaining positions conflict-free colors every interval with $\lceil\log_2|X|\rceil$ colors and give the stated bound.
\end{proof}

\begin{theorem}\label{thm:split-cfon}
Let $G=(C\dot\cup I,E)$ be connected, clique-side circular-convex, $|C|\ge3$. Put $\mathcal S=\{N(x):x\in I\}$. Then
$
 \chi_{\ON}^*(G)\le a^*(C,\mathcal S),\qquad
 \chi_{\ON}(G)\le a(C,\mathcal S)\le1+\lceil\log_2|C|\rceil.
$
These are polynomial-time constructive upper bounds when the cyclic order is supplied.
\end{theorem}
\begin{proof}
For the partial variant, we apply Lemma~\ref{lem:anchor} on $C$ and leave $I$ uncolored. Every independent vertex sees its interval witness. Every clique vertex sees one of the two globally singleton clique colors at a distinct clique vertex.

For the full variant, we choose the anchored full coloring on $C$, with anchor $p$ carrying its new color and $q\in C\setminus\{p\}$ carrying a color $\beta$ unique on $C\setminus\{p\}$. The coloring of $C\setminus\{p\}$ uses at least two colors: the hyperedge $C\setminus\{p\}$ cannot be conflict-free monochromatic. Choose a color $\gamma\ne\beta$ occurring there and assign $\gamma$ to every vertex in $I$. Every independent vertex retains its clique witness. Each clique vertex except $p$ sees the anchor color once, while $p$ sees $\beta$ once because no vertex of $I$ has color $\beta$.
\end{proof}

\begin{theorem}\label{thm:starlike-cfon}
Suppose $G$ is starlike with $|C|\ge3$, and $S_1,\ldots,S_m$ are circular intervals in one cyclic ordering of $C$. Then, writing $\mathcal S=\{S_i:i\in[m]\}$,
\[
 \chi_{\ON}^*(G)\le a^*(C,\mathcal S),\qquad
 \chi_{\ON}(G)\le a(C,\mathcal S)\le1+\lceil\log_2|C|\rceil.
\]
\end{theorem}
\begin{proof}
For partial CFON, we color only $C$ as in Lemma~\ref{lem:anchor}. A vertex in $Q_i$ has open neighborhood $(Q_i\setminus\{v\})\cup S_i$, and the peripheral vertices are uncolored, so it sees its witness in $S_i$.

For full CFON, anchor $p$ and let $\beta$ be the unique color at $q\in C\setminus\{p\}$. For each $Q_i$, if $p\in S_i$, we choose a color occurring in $C\setminus\{p\}$ different from $\beta$ and color $Q_i$ monochromatically with it. If $p\notin S_i$, choose a color $\beta_i$ unique on $S_i$ and color $Q_i$ monochromatically with another color from $C\setminus\{p\}$. In either case a peripheral vertex retains a clique-side unique witness, even though its own peripheral clique has repeated colors. Every central vertex other than $p$ sees the new anchor color exactly once, and $p$ sees $\beta$ exactly once: no peripheral clique adjacent to $p$ receives $\beta$.
\end{proof}

\begin{proposition}\label{prop:tight-cfon}
For each $n\ge3$ there is a clique-side circular-convex split graph $G_n$, also starlike, with $|C|=n$ and
\(
 \chi_{\ON}^*(G_n)=\chi_{\ON}(G_n)=1+\lceil\log_2n\rceil.
\)
\end{proposition}
\begin{proof}
We add one independent vertex with neighborhood $S$ for every nonempty circular interval $S\subseteq C$. The singleton intervals force all clique vertices to receive nonzero colors even in a partial coloring. The independent vertex with neighborhood $C$ forces a color occurring uniquely on $C$, at some $p$. On the remaining $n-1$ linearly ordered positions, every nonempty interval is represented, and the color of $p$ appears nowhere. A conflict-free coloring of all nonempty intervals on $n-1$ ordered points requires $\lceil\log_2 n\rceil$ colors: recursively removing the unique witness of the whole interval shows that $k$ colors cover at most $2^k-1$ points. Then Theorem~\ref{thm:split-cfon} gives the matching upper bound.
\end{proof}

\subsection{Bounds on WPC graphs}
For every bag $B$, we define $\mathcal S_B=\{S_{B,D}:D\in N_T(B)\}$ and let $a_B^*=a^*(B,\mathcal S_B)$, $a_B=a(B,\mathcal S_B)$, with the singleton-bag convention from Lemma~\ref{lem:anchor}. Root $T$ at a nonsingleton bag if one exists, otherwise at a leaf. Let $d(B)$ be its depth and put
\[
 L_j^*=\max_{d(B)\equiv j\ (\mathrm{mod}\,3)}a_B^*,\qquad
 L_j=\max_{d(B)\equiv j\ (\mathrm{mod}\,3)}a_B \qquad(j=0,1,2),
\]
with empty maxima zero.

\begin{theorem}\label{thm:wpc-cfon}
For a connected bagwise circular-convex WPC graph with a supplied partition tree and cyclic orders,
\[
 \chi_{\ON}^*(G)\le\sum_{j=0}^2L_j^*,\qquad
 \chi_{\ON}(G)\le\min\left\{\sum_{j=0}^2L_j,\ 1+\sum_{j=0}^2L_j^*\right\}.
\]
Consequently, if $b=\max_{B\in\mathcal B}|B|$, then
\(
 \chi_{\ON}^*(G)\le\chi_{\ON}(G)\le3\bigl(1+\lceil\log_2 b\rceil\bigr).
\)
The bound is computable in polynomial time from the supplied representation.
\end{theorem}
\begin{proof}
We use pairwise disjoint colors for bags at depths congruent to $0,1,2$ modulo three, reusing the colors among bags in the same residue class. On nonsingleton bags we apply Lemma~\ref{lem:anchor} and on singleton bags we  color the sole vertex with a positive color of that depth class. A nonsingleton bag vertex sees a unique color of its own bag on another bag vertex, and adjacent bags use different colors. If a bag $B=\{v\}$ is singleton and not the root, its parent $D$ has a boundary $S_{D,B}$ containing a unique color and all children of $B$ use the depth class $d(B)+1$ instead of the parent's class $d(B)-1$, so the parent's witness is not repeated. If every bag is singleton, the root was chosen as a leaf and its unique neighbor supplies a witness. This proves the partial bound; the same argument with full local colorings proves $\sum_j L_j$. Alternatively, we recolor all uncolored vertices with one fresh color. The logarithmic estimate follows from Lemma~\ref{lem:anchor}.
\end{proof}
% \begin{remark}
% The quantities $a_B,a_B^*$ are \emph{upper-bound parameters}; Theorem~\ref{thm:wpc-cfon} is not an exact optimization algorithm for CFON on arbitrary bagwise circular-convex WPC graphs. The lower-bound family in Proposition~\ref{prop:tight-cfon} is itself a star-partition WPC graph and shows the order $\Omega(\log b)$ is unavoidable.
% \end{remark}
\section{CFCN and CFCN*: bounds, characterizations, and hardness}\label{sec:cfcn}
\begin{theorem}\label{thm:starlike-cfcn}
Let $G$ be a connected starlike graph with central clique $C$ and indexed peripheral neighborhoods $S_1,\ldots,S_m$. Then $\chi_{\CN}^*(G)\le2$. Moreover,
\[
 \chi_{\CN}^*(G)=1\quad\Longleftrightarrow\quad
 \text{$G$ has a universal vertex, or each $c\in C$ belongs to exactly one $S_i$.}
\]
Thus the exact partial parameter can be computed in polynomial time from a starlike partition. For a connected noncomplete starlike graph with an edge, $2\le\chi_{\CN}(G)\le3$.
\end{theorem}
\begin{proof}
First, we choose $p\in C$ and color it 1; for each $i$ with $p\notin S_i$, we color one arbitrary vertex in $Q_i$ with 2; leave all other vertices uncolored. Each central vertex sees color 1 exactly once. A vertex of $Q_i$ sees color 1 at $p$ if $p\in S_i$, and otherwise sees color 2 exactly once within its own clique. This gives the upper bound.

For the characterization, suppose $D$ is an efficient dominating set. Since $C$ is a clique, $|D\cap C|\le1$. If $D\cap C=\{p\}$ and $p\notin S_i$ for some $i$, a vertex of $Q_i$ must be chosen into $D$; then any $u\in S_i$ sees both $p$ and the selected peripheral vertex, a contradiction. Hence $p$ is adjacent to all peripheral cliques and is universal. If $D\cap C=\varnothing$, exactly one vertex from each $Q_i$ lies in $D$, and a central vertex is dominated exactly once if and only if it belongs to exactly one indexed set $S_i$. Conversely, a universal vertex alone, or one representative from each $Q_i$ when the $S_i$ partition $C$ with multiplicities, is an efficient dominating set.

Coloring all previously uncolored vertices with color 3 proves the full upper bound; a noncomplete graph with an edge cannot have full CFCN number one.
\end{proof}
\begin{theorem}
\label{thm:starlike-full-cfcn}
Let $G$ be a connected starlike graph with central clique $C$ and
pairwise anticomplete peripheral cliques $Q_1,\ldots,Q_m$.
Every vertex of $Q_i$ has the same nonempty neighborhood
$S_i\subseteq C$. The sets $S_1,\ldots,S_m$ are indexed, so
repetitions are counted separately. Then
 $\chi_{\CN}(G)=
 \begin{cases}
 1,& |V(G)|=1,\\
 2,& |V(G)|\ge2\text{ and at least one of (a)--(c) holds},\\
 3,& \text{otherwise},
 \end{cases}$
where the conditions are:
\begin{enumerate}[label=\textup{(\alph*)},leftmargin=2.3em]
 \item $G$ has a universal vertex; equivalently, some $p\in C$
       belongs to every $S_i$;
 \item the indexed family $(S_1,\ldots,S_m)$ partitions $C$,
       that is,
       $ \sum_{i=1}^m \mathbf{1}_{\{c\in S_i\}}=1
        \quad\text{for all }c\in C;
       $
 \item $|C|=2$ and every $S_i$ is a singleton.
\end{enumerate}
Consequently, full CFCN is solvable in polynomial time on starlike
graphs when a starlike partition is supplied. In particular, if
$|C|\ge3$, condition \textup{(c)} is unnecessary.
\end{theorem}
\begin{proof}
The case $|V(G)|=1$ is immediate. If $|V(G)|\ge2$, one color
cannot work: $G$ is connected and has an edge, and the closed
neighborhood of an endpoint of that edge contains at least two
vertices. Theorem~\ref{thm:starlike-cfcn} gives
$\chi_{\CN}(G)\le3$. We therefore characterize when two colors
suffice.

\emph{Sufficiency.}
Suppose \textup{(a)} holds and let $p$ be universal. Color $p$
with color 1 and all other vertices with color 2. Every closed
neighborhood contains $p$ and therefore has exactly one occurrence
of color 1.

Suppose \textup{(b)} holds. Color every vertex of $C$ with
color 1. In each $Q_i$, choose one vertex and color it 2; color
all other vertices of $Q_i$ with 1. Every vertex of $Q_i$
sees exactly one color-2 vertex in $N[v]$, namely the chosen
vertex of its own clique. Every central vertex $c$ belongs to
exactly one $S_i$, so it also sees exactly one color-2 vertex.

Finally, suppose \textup{(c)} holds and write $C=\{a,b\}$.
Color $a$ with 1 and $b$ with 2. Color each $Q_i$
monochromatically with color 2 when $S_i=\{a\}$, and with
color 1 when $S_i=\{b\}$. The vertex $a$ sees color 1
exactly once in $N[a]$, at $a$, and $b$ sees color 2 exactly
once, at $b$. A vertex in a clique with boundary $\{a\}$
sees color 1 exactly once, at $a$; the symmetric statement
holds for boundary $\{b\}$. Thus this coloring is valid.

\emph{Necessity when $|C|\ge3$.}
Let $\gamma$ be a full CFCN coloring with colors 1 and 2.
If both colors occur at least twice on $C$, every central
vertex sees both colors at least twice in its closed
neighborhood, a contradiction. Thus either $C$ is
monochromatic, or one color occurs on exactly one central
vertex.

First suppose that all of $C$ has color 1. Since $|C|\ge3$,
each central vertex must see color 2 exactly once among its
peripheral neighbors. Moreover, for every $i$, at least one
vertex of $Q_i$ must have color 2: otherwise
$N[v]=Q_i\cup S_i$ is monochromatic and contains at least
two vertices for each $v\in Q_i$. Pick any $c\in S_i$.
As $c$ must see exactly one color-2 peripheral vertex,
$Q_i$ has \emph{exactly} one vertex of color 2.
Consequently, the number of color-2 vertices in the
closed neighborhood of a central vertex $c$ equals
$|\{i:c\in S_i\}|$. It must be one, proving \textup{(b)}.

Next suppose that one central vertex $p$ has color 2 and
all other central vertices have color 1. Since at least
two central vertices have color 1, every central vertex
can only use color 2 as a unique witness. This color already
appears exactly once on $C$, at $p$. Thus no peripheral
vertex can have color 2, because each $S_i$ is nonempty.
If some $S_i$ omitted $p$, every vertex of $Q_i$ would have
a monochromatic closed neighborhood of size at least two.
Hence $p\in S_i$ for every $i$, which is \textup{(a)}.

\emph{Necessity when $|C|=2$.}
Write $C=\{a,b\}$. Every nonempty $S_i$ is one of
$\{a\}$, $\{b\}$, and $\{a,b\}$. If at least one type is
absent, then either $a$ or $b$ is universal (when type
$\{b\}$ or $\{a\}$ is absent, respectively), or all
$S_i$ are singletons (when type $\{a,b\}$ is absent).
Thus it remains to show that no valid two-coloring exists
when all three types occur.

Suppose first that $a,b$ have the same color, say 1.
Every peripheral clique with boundary $\{a\}$ or
$\{a,b\}$ must contain at least one color-2 vertex;
otherwise its vertices have monochromatic closed
neighborhoods of size at least two. Both cliques attach
to $a$, so $N[a]$ contains at least two vertices of
color 2, while $a$ and $b$ already contribute two
vertices of color 1. Thus $a$ has no uniquely occurring
color, a contradiction.

Now suppose $\gamma(a)=1$ and $\gamma(b)=2$. A peripheral
clique with boundary $\{a\}$ must contain a color-2
vertex; otherwise all colors in each of its closed
neighborhoods are 1. Therefore, color 2 occurs at
least twice in $N[a]$ (at $b$ and this peripheral
vertex). For $a$ to be satisfied, color 1 must appear
exactly once in $N[a]$, forcing every peripheral vertex
adjacent to $a$ to have color 2. In particular, all
vertices in every clique of boundary $\{a,b\}$ have
color 2. Symmetrically, a clique with boundary $\{b\}$
forces every peripheral vertex adjacent to $b$ to
have color 1. This is impossible for any clique
with boundary $\{a,b\}$. The contradiction finishes
the proof.

Finally, conditions \textup{(a)}--\textup{(c)} can be checked
by scanning the indexed neighborhoods $S_i$, so the
classification yields a polynomial-time algorithm.
\end{proof}

% \begin{remark}
% \label{rem:starlike-small-counterexample}
% The frequently proposed condition $|C|\le2$ is
% \emph{not} sufficient for full two-color CFCN. For example,
% let $C=\{a,b\}$ and let three singleton peripheral
% cliques have neighborhoods $\{a\}$, $\{b\}$, and
% $\{a,b\}$, respectively. All three boundary types
% occur, so Theorem~\ref{thm:starlike-full-cfcn} gives
% $\chi_{\CN}(G)=3$, even though $|C|=2$.
% \end{remark}

\subsection{WPC graphs: sharp bounds and hardness}
\begin{theorem}\label{thm:wpc-cfcn}
Every nonempty WPC graph has
\(
 \chi_{\CN}^*(G)\le2,\;\chi_{\CN}(G)\le3.
\)
Both colorings are computable in polynomial time from a supplied partition forest. The bounds are simultaneously tight on the 3-sun.
\end{theorem}
\begin{proof}
We first bipartition the partition tree into even- and odd-depth bags. Then we choose one representative $r_B$ per bag, give it color 1 or 2 according to parity, and leave all other vertices uncolored. For every $v\in B$, the representative $r_B$ is in $N[v]$ and no other bag of the same parity is adjacent to $B$, so its color occurs exactly once in $N[v]$. We give all uncolored vertices a fresh color 3 for the full version.

For sharpness, the 3-sun has a central clique $\{c_1,c_2,c_3\}$ and independent vertices $x_i$ adjacent to $c_i,c_{i+1}$ (indices modulo 3). The indexed neighborhoods do not partition the central clique and no vertex is universal, so the partial optimum is 2 by Theorem~\ref{thm:starlike-cfcn}. A full three-color CFCN coloring gives $c_1$ color 1, $c_2,c_3$ color 2, and every $x_i$ color 3. For two colors, the central clique cannot contain at least two of each color; if it is monochromatic, satisfying each central vertex forces adjacent independent vertices to have opposite colors around a 3-cycle, impossible. If the central clique has color multiplicities $1,2$, a peripheral vertex adjacent to the uniquely colored clique vertex sees both colors at least twice in its closed neighborhood, again impossible. Hence the full optimum is 3.
\end{proof}

\begin{theorem}\label{thm:wpc-cfcnstar-hard}
Deciding $\chi_{\CN}^*(G)=1$ is NP-complete on WPC graphs even when the partition tree is a supplied star and every peripheral bag has two vertices.
\end{theorem}
\begin{proof}
Given an X3C instance $U=\{u_1,\dots,u_{3q}\}$ and triples $S_i\subseteq U$ with every element in at least one set, let $C=U$ be a clique. For every triple $S_i$, introduce a clique $Q_i=\{a_i,b_i\}$ with $a_i$ adjacent to exactly $S_i$ in $C$, $b_i$ having no central neighbor, and no edges between different $Q_i$. The partition tree is a star centered at $C$, with complete bipartite boundary edges $S_i\times\{a_i\}$, hence WPC.

Suppose $D$ is efficient dominating. Since $N[b_i]=\{a_i,b_i\}$, exactly one of $a_i,b_i$ is selected. No central vertex $u$ can be in $D$: choose $S_i\ni u$ and note that $N[a_i]$ would then contain $u$ and the selected vertex from $\{a_i,b_i\}$. Hence $D\cap C=\varnothing$. Each central $u$ must see exactly one selected $a_i$ with $u\in S_i$, so selected $a_i$ represent an exact cover. Conversely, select $a_i$ for sets in an exact cover and $b_i$ for other sets; the resulting $D$ is efficient dominating. Membership in NP is immediate.
\end{proof}

\begin{theorem}\label{thm:wpc-cfcnfull-hard}
Deciding whether $\chi_{\CN}(G)\le2$ is NP-complete on WPC graphs. Consequently, determining the full conflict-free
closed-neighborhood chromatic number of a chordal
graph is NP-hard.
\end{theorem}
\begin{proof}
We use the same X3C graph. Given an exact cover $\mathcal S'$, color all central vertices 1 and put $(c(a_i),c(b_i))=(2,1)$ for $S_i\in\mathcal S'$ and $(1,2)$ otherwise. Each central $u$ sees exactly one color-2 neighbor, while $a_i$ and $b_i$ each see color 2 exactly once in their closed neighborhoods.

Conversely, suppose a full two-color CFCN coloring exists. Because $N[b_i]=\{a_i,b_i\}$, the vertices $a_i,b_i$ have opposite colors. In the central clique $C$, neither color can occur at least twice, because then no central vertex would have a unique closed-neighborhood color. If the clique is not monochromatic, its colors have multiplicities $1,|C|-1$ (recall $|C|\ge3$). Take a triple containing the unique minority-colored central vertex. The corresponding $N[a_i]$ contains that minority once and the majority twice from $C$, plus the opposite-colored pair $a_i,b_i$, giving multiplicities 2 and 3; this fails the conflict-free condition. Thus $C$ is monochromatic; assume it has color 1. Each central vertex then has at least three occurrences of color 1 in its closed neighborhood, so it must see exactly one vertex $a_i$ of color 2 among its incident triples. The $S_i$ for which $a_i$ has color 2 form an exact cover. Membership in NP is immediate.
\end{proof}
% \begin{remark}
% The two hardness theorems above are proved \emph{separately}: hardness of one-color partial CFCN alone cannot be used to assert hardness of two-color full CFCN. The X3C triples are arbitrary and need not be circular intervals, so these reductions do not establish hardness on bagwise circular-convex WPC graphs.
% \end{remark}
\section{Rainbow vertex connection on split, starlike and well partitioned chordal graphs}\label{sec:rainbow}
We use the rainbow vertex-connection conventions fixed in the preliminaries; endpoints of the paths impose no color restriction. 

We first establish a basic lemma regarding $\rvc$ and $\srvc$ of split graphs.

\begin{lemma}\label{lem:split-rainbow}
For every connected split graph $H$,
$\rvc(H)=\srvc(H)$.
\end{lemma}
\begin{proof}
Every pair at distance at most two has a rainbow shortest path under
any vertex coloring.  For a pair $x,y$ at distance three, both are
independent vertices with disjoint neighborhoods in the clique.
On a rainbow $x$--$y$ path, its first and last internal clique
vertices have different colors.  The edges between those two clique
vertices give a rainbow geodesic $x-a-b-y$.
\end{proof}

\begin{theorem}\label{lem:local-rainbow}
Let $H=(C\dot\cup I,E)$ be a connected circular-convex split graph with respect to $C$, together with $|I|\ge2$. Let
\(
 P=\{p\in C:N(x)=\{p\}\text{ for some }x\in I\},\quad s=|P|,
\)
and let $\delta=1$ if two neighborhoods $N(x),N(y)$, $x,y\in I$, are disjoint, and $\delta=0$ otherwise. Then
\(\rvc(H)=\srvc(H)=\max\{s,1+\delta\}.\)
The optimum is polynomial-time computable from the circular representation.
\end{theorem}

\begin{proof}
We first establish the lower bound. Since $|I|\ge2$,
the graph is not complete, so $\rvc(H)\ge1$.
If two vertices of $I$ have disjoint neighborhoods,
their distance is three, which requires at least
two colors. Moreover, for any distinct $p,q\in P$,
there exist $x,y\in I$ with $N(x)=\{p\}$ and
$N(y)=\{q\}$. Every $x$--$y$ path contains
$p$ and $q$ as internal vertices, so they must
receive different colors. Hence,
$\rvc(H)\ge\max\{s,1+\delta\}$.

We now prove the upper bound. If $\delta=0$,
every two vertices of $I$ have a common neighbor.
Thus $\diam(H)\le2$, and one color is sufficient.
This settles the case $s\le1$ and $\delta=0$.

Suppose $\delta=1$. We construct a coloring
of $C$ as follows.

\noindent\textbf{Case 1: $s\ge2$.}
Give the vertices of $P$ distinct colors
$1,\ldots,s$. In the circular ordering of $C$,
color the remaining vertices between every two
consecutive vertices of $P$ alternately with
colors $1$ and $2$.

If $x\in I$ satisfies $N(x)\cap P=\varnothing$,
then $N(x)$ is a circular interval lying between
two consecutive vertices of $P$. Since
$|N(x)|\ge2$, it contains both colors $1$ and $2$.

\noindent\textbf{Case 2: $s=1$.}
Let $P=\{p\}$. Give $p$ color $1$ and color the
remaining vertices of $C$ alternately with colors
$1$ and $2$ in the linear order obtained by
removing $p$. Again, every neighborhood avoiding
$P$ contains both colors $1$ and $2$.

In both cases, consider $x,y\in I$ with
$N(x)\cap N(y)=\varnothing$. If both neighborhoods
contain vertices of $P$, choose one such vertex
from each; they have different colors. Otherwise,
at least one neighborhood avoids $P$ and hence
contains both colors $1$ and $2$. We can therefore
choose $a\in N(x)$ and $b\in N(y)$ with
different colors. Since $C$ is a clique,
$x-a-b-y$ is a rainbow shortest path.

\noindent\textbf{Case 3: $s=0$.}
Color the vertices of $C$ alternately with
colors $1$ and $2$ in their circular order.
If $|C|$ is even, every circular interval of
size at least two contains both colors.
If $|C|$ is odd, only the pair at which the
coloring wraps around can be monochromatic.
Since every vertex of $I$ has at least two
neighbors, two disjoint neighborhoods cannot
both be monochromatic. Hence, for any
$x,y\in I$ with disjoint neighborhoods,
we can choose differently colored vertices
$a\in N(x)$ and $b\in N(y)$, giving the
rainbow shortest path $x-a-b-y$.

In all three cases, every pair of independent
vertices with disjoint neighborhoods has a
rainbow shortest path. All other vertex pairs
have distance at most two and hence automatically
have rainbow shortest paths.
Thus the constructed coloring uses at most
$\max\{s,1+\delta\}$ colors, giving
$\srvc(H)\le\max\{s,1+\delta\}.$
Combining this with the lower bound and
$\rvc(H)\le\srvc(H)$ yields
$\rvc(H)=\srvc(H)=\max\{s,1+\delta\}.$
Finally, $P$ and $\delta$ can be determined by
examining the neighborhoods of vertices in $I$.
The coloring can also be constructed in polynomial
time from the given circular ordering.
\end{proof}

\begin{theorem}\label{cor:rainbow-star}
Let $G$ be connected and noncomplete, with a starlike decomposition consisting of central clique $C$ and peripheral cliques $Q_1,\dots,Q_m$ having circular-interval central neighborhoods $S_i$. Let $s$ count distinct vertices $p\in C$ such that some $S_i=\{p\}$. Then
\(
 \rvc(G)=\srvc(G)=\max\{\diam(G)-1,s\}.
\)
\end{theorem}
\begin{proof}
Every pair of vertices at distance three belongs to distinct peripheral cliques with disjoint central neighborhoods, and its shortest path has the form $x-a-b-y$, where $a\in S_i$, $b\in S_j$. The coloring argument of Theorem~\ref{lem:local-rainbow}, applied to the circular intervals $S_i$, guarantees differently colored choices $a,b$. All other pairs have distance at most two. For the lower bound, distinct singleton-support vertices lie internally on paths between vertices of their associated peripheral cliques, forcing distinct colors, and the diameter lower bound applies. Consequently the split-component proof extends verbatim to the starlike construction, even when the peripheral cliques have more than one vertex.
\end{proof}
\subsection{RVC and SRVC Bounds on WPC Graphs}
Before presenting the results on WPC graphs, we define few things.

% For a leaf bag $L$ with neighbor $P(L)$, call it \emph{nontrivial} if $S_{L,P(L)}\subsetneq L$ and set
%  $
%  \mathcal L_{\rm nt}=\{L:\deg_T(L)=1,\ S_{L,P(L)}\subsetneq L\}.
% $
% Define
% $
%  K(G,T)=\sum_{\substack{B\in V(T)\\\deg_T(B)\ge2}}\rvc(H_B).
% $ Also, we partition nontrivial leaves into
% $
%  \mathcal L_1=\{L\in\mathcal L_{\rm nt}:|S_{L,P(L)}|=1\},\qquad
%  \mathcal L_{\ge2}=\{L\in\mathcal L_{\rm nt}:|S_{L,P(L)}|\ge2\}.
% $
% Put $s=|\mathcal L_1|$ and $f=|\mathcal L_{\ge2}|$.
% Requires the existing manuscript definitions of \rvc, \srvc, WPC partition trees,
% attachment boundaries S_{B,D}, and split components H_B.
% Also uses the previously established split-graph equality
% \rvc(H)=\srvc(H), cited as Lemma~\ref{lem:split-rainbow}.

Let $G$ be a connected WPC graph with a supplied partition tree $T$.
For adjacent bags $B,D$, write $S_{B,D}\subseteq B$ for the attachment
boundary in $B$. For a nonleaf bag $B$, let $H_B$ be the split graph
whose clique is $B$ and which has one independent marker for each
$D\in N_T(B)$, adjacent precisely to $S_{B,D}$.
We use the convention $\rvc(G)=\srvc(G)=0$ when $G$ is complete.

A leaf bag $L$ with neighbor $P(L)$ is \emph{nontrivial} if
$S_{L,P(L)}\subsetneq L$. Define
\(
\mathcal L_1=\{L:\deg_T(L)=1,\ S_{L,P(L)}\subsetneq L,
                     \ |S_{L,P(L)}|=1\},
\)
and
\(
\mathcal L_{\ge2}=\{L:\deg_T(L)=1,\ S_{L,P(L)}\subsetneq L,
                     \ |S_{L,P(L)}|\ge2\}.
\)
Put $s=|\mathcal L_1|$, $f=|\mathcal L_{\ge2}|$, and
\(
K(G,T)=\sum_{\deg_T(B)\ge2}\rvc(H_B),\qquad
\Psi(G,T)=K(G,T)+s+\min\{2,f\}.
\)
The earlier one-new-color-per-nontrivial-leaf estimate is
\(
\Phi(G,T)=K(G,T)+s+f.
\)
The bound below improves this to $\Psi(G,T)$, and branching permits a
further saving.

We call the subtree induced by the nonleaf bags the \emph{nonleaf subtree}
$T^\circ$. If there are no nonleaf bags, $T^\circ$ is empty.
We say that $T$ has \emph{branching nonleaf bags} if $T^\circ$
contains a vertex of degree at least three; equivalently, no path
of $T$ contains all nonleaf bags.

\begin{theorem}\label{thm:rainbow-improve}
Every connected WPC graph satisfies
\(\rvc(G)\le\srvc(G)\le\Psi(G,T):=K(G,T)+s+\min\{2,f\}\).
% More precisely, $\Phi(G,T)-\Psi(G,T)=\max\{0,f-2\}$.
\end{theorem}
\begin{proof}
For each nonleaf bag $B$, we use a strong rainbow vertex
coloring of its split component $H_B$ with
$\rvc(H_B)$ colors. By
Theorem~\ref{lem:split-rainbow}, such a coloring
exists. We use different sets of colors for
different nonleaf bags. This requires
$K(G,T)$ colors.

For each leaf bag in $\mathcal L_1$, we assign
a distinct new color to its unique attachment
vertex.

Now consider the leaf bags in $\mathcal L_{\ge2}$.
If there is only one such leaf bag, we assign one
new color to all its attachment vertices.
If there are at least two, we introduce two new
colors $\alpha$ and $\beta$. In each such leaf
bag, color one attachment vertex with $\alpha$
and another with $\beta$. The remaining
attachment vertices may receive either color.
The same two colors are used in all these
leaf bags.

We now show that the coloring is a strong
rainbow vertex coloring of $G$.

Consider any two vertices $x,y$. A shortest
$x$--$y$ path follows the unique path between
their bags in the partition tree. Within each
nonleaf bag, the coloring of its split component
ensures that the required internal vertices
have different colors. Since different nonleaf
bags use different colors, no color is repeated
between them. Furthermore, only the leaf bags containing
$x$ and $y$ can contribute internal vertices.
Thus, at most two leaf bags contribute.
If both have at least two attachment vertices,
we choose a vertex colored $\alpha$ in one
leaf bag and a vertex colored $\beta$ in the
other. Their colors are different.
These choices give a shortest path because
every bag is a clique and the edges between
adjacent bags are complete between their
attachment boundaries.
Hence every pair of vertices has a rainbow
shortest path. The total number of colors
used is at most
$K(G,T)+s+\min\{2,f\}.$
Therefore,
$\rvc(G)\le\srvc(G)\le\Psi(G,T).$
\end{proof}

The preceding bound can be sharpened whenever $f>0$ and the
nonleaf subtree branches.

\begin{theorem}
\label{thm:rainbow-branch-improve}
Let $G$ be a connected WPC graph with $f\ge1$.
If $T^\circ$ contains a branching vertex, then
\(
\rvc(G)\le\srvc(G)\le\Psi(G,T)-1.
\)
\end{theorem}

\begin{proof}
We begin with the coloring constructed in
Theorem~\ref{thm:rainbow-improve}.
Recall that different nonleaf bags use different
sets of colors, and these colors are also
different from the colors assigned to leaf bags.

Since $T^\circ$ contains a branching vertex,
there is a nonleaf bag $B$ such that at least
three components of $T-B$ contain nonleaf bags.
Choose three such components $T_1,T_2,T_3$,
and a nonleaf bag $A_i$ in each $T_i$.
Let $\gamma_i$ be a color used on a vertex
of $A_i$, for $i=1,2,3$.

The colors $\gamma_1,\gamma_2,\gamma_3$
are distinct and do not occur in any other
nonleaf bag. Moreover, a path in $T$ can
visit at most two of the three components
$T_1,T_2,T_3$.

We show that one of the fresh colors
assigned to the leaf bags in
$\mathcal L_{\ge2}$ can be eliminated.

\medskip
\noindent
\textbf{Case 1: $f=1$.}

Let $L$ be the only leaf bag in
$\mathcal L_{\ge2}$. At least two of
$T_1,T_2,T_3$ do not contain $L$.
Suppose these components are $T_i$ and $T_j$.

Instead of using a fresh color on the
attachment boundary of $L$, color one
attachment vertex with $\gamma_i$ and
another with $\gamma_j$. The remaining
attachment vertices may receive either color.

Consider a shortest path from a vertex of
$L$ to a vertex outside $L$. Its route in
$T$ cannot contain both $A_i$ and $A_j$.
Hence at least one of these source bags
lies outside the route.

We choose the attachment vertex of $L$
whose color comes from that source bag.
Since this color does not occur in any
nonleaf bag on the route, it cannot repeat
an internal color of the path. It is also
different from all other leaf colors.

As the edges between adjacent bags are
complete between their attachment
boundaries, this choice does not increase
the length of the path. Thus a rainbow
shortest path still exists.

Paths whose endpoints are outside $L$
remain unchanged. Therefore, the fresh
color previously used for $L$ is no
longer needed.

\medskip
\noindent
\textbf{Case 2: $f\ge2$.}

In the original coloring, two fresh
colors $\alpha$ and $\beta$ are used on
the attachment boundaries of all leaf
bags in $\mathcal L_{\ge2}$.

We retain $\beta$ and replace $\alpha$
with colors already used in nonleaf bags.

For each $L\in\mathcal L_{\ge2}$, we define
\(
\gamma(L)=
\begin{cases}
\gamma_2,&L\in T_1,\\
\gamma_3,&L\in T_2,\\
\gamma_1,&\text{otherwise}.
\end{cases}
\)

We color one attachment vertex of $L$
with $\gamma(L)$ and another with $\beta$.
All remaining attachment vertices may
receive either of these two colors.

We claim that, for any two distinct leaf
bags $L,R\in\mathcal L_{\ge2}$, at least
one of the source bags of $\gamma(L)$ and
$\gamma(R)$ lies outside their connecting
route in $T$.

To see this, first suppose that $L$ and
$R$ lie in the same component of $T-B$.
Their connecting route stays inside that
component, whereas their borrowed colors
come from another component. Thus both
borrowed colors are safe.

Now suppose $L$ and $R$ lie in different
components of $T-B$. Their connecting
route passes through $B$ and visits only
the two components containing $L$ and $R$.
The replacement rule assigns colors from
$T_1,T_2,T_3$ in the cyclic order
\(
T_1\to T_2\to T_3\to T_1.
\)
Therefore, the two leaves cannot both
borrow colors from each other's component.
Hence at least one borrowed color comes
from a source bag outside their route.
The same argument applies when one leaf
is adjacent directly to $B$ or lies in
another component of $T-B$.

We now construct a rainbow shortest path
between vertices of $L$ and $R$.
Choose an attachment vertex in one leaf whose color is taken from a nonleaf bag not lying on the path between \(L\) and \(R\). In the other leaf, choose an attachment vertex colored \(\beta\). Since the chosen nonleaf bag does not lie on this path, its color cannot appear on any nonleaf vertex of the path. Also, \(\beta\) is not used in any nonleaf bag. Therefore, the two chosen attachment vertices do not create any repeated color, and the resulting shortest path is rainbow.

If only one endpoint belongs to a leaf
bag in $\mathcal L_{\ge2}$, we choose its
attachment vertex colored $\beta$.
All paths with neither endpoint in this
family remain unchanged.
Hence the recoloring preserves strong
rainbow vertex connectivity.

In both cases, we eliminate one fresh
color without introducing any new color.
Therefore, the resulting strong rainbow
vertex coloring uses at most
\(
\Psi(G,T)-1
\)
colors. Consequently,
\(
\rvc(G)\le\srvc(G)\le\Psi(G,T)-1.
\)
\end{proof}

Define
\(
\eta(T,f)=
\begin{cases}
1,&f\ge1\text{ and }T^\circ\text{ has a branching vertex},\\
0,&\text{otherwise},
\end{cases}
\qquad
\widehat\Psi(G,T)=\Psi(G,T)-\eta(T,f).
\)

\begin{corollary}\label{cor:wpc-unified}
Every connected WPC graph satisfies
\(
\rvc(G)\le\srvc(G)\le\widehat\Psi(G,T)\le\Psi(G,T).
\)
\end{corollary}

Let $\ell(G)$ be the number of cut vertices, and set
$D(G)=\max\{\diam(G)-1,\ell(G)\}$. The standard lower bounds
and Corollary~\ref{cor:wpc-unified} yield
\(
D(G)\le\rvc(G)\le\srvc(G)\le\widehat\Psi(G,T).
\)
\begin{theorem}\label{thm:psi-certificate}
If $D(G)=\widehat\Psi(G,T)$, then
\(
\rvc(G)=\srvc(G)=\widehat\Psi(G,T).
\)
\end{theorem}

\begin{proof}
The diameter and cut-vertex lower bounds give
$D(G)\le\rvc(G)$, while Corollary~\ref{cor:wpc-unified}
gives $\srvc(G)\le\widehat\Psi(G,T)$. Hence
\(
D(G)\le\rvc(G)\le\srvc(G)\le\widehat\Psi(G,T).
\)
If the outer terms agree, all inequalities are equalities.
\end{proof}

We now give structural exactness classes according to the value of
$f$ and the shape of $T^\circ$.

\subsection{Exactness Classes for the Bound}

\subsubsection{Case I: $f=0$ (arbitrary branching is allowed)}

For each nonleaf bag $B$ define
\(
P_B=\{v\in B:S_{B,D}=\{v\}\text{ for some }D\in N_T(B)\},
\qquad t_B=|P_B|.
\)

\begin{lemma}\label{lem:wpc-cut-count}
If $T$ has at least two bags, then
\(
\ell(G)=s+\sum_{\deg_T(B)\ge2}t_B.
\)
\end{lemma}

\begin{proof}
A vertex $v$ of a nonleaf bag $B$ is a cut vertex if and only if
$S_{B,D}=\{v\}$ for some neighboring bag $D$: deleting such a
vertex disconnects the part of the graph beyond $D$, whereas
if no attachment boundary is the singleton $\{v\}$, every
boundary remains nonempty after deleting $v$ and the remaining
graph is connected. In a leaf bag $L$, its attachment vertex
is a cut vertex precisely when its attachment boundary is a
proper singleton. Such leaves are counted by $s$.
Summing over all bags proves the formula.
\end{proof}

\begin{theorem}
\label{thm:wpc-cut-exact}
Suppose $G$ is connected, $T$ has at least two bags,
$f=0$, and every nonleaf bag $B$ satisfies
$\rvc(H_B)=t_B$. Then
\(
\rvc(G)=\srvc(G)=\Psi(G,T)=\ell(G).
\)
In particular, the nonleaf subtree may branch arbitrarily.
\end{theorem}

\begin{proof}
Since $f=0$ and $\rvc(H_B)=t_B$ for every nonleaf bag,
\(
\Psi(G,T)=K(G,T)+s
=\sum_{\deg_T(B)\ge2}t_B+s
=\ell(G),
\)
where the last equality follows from Lemma~\ref{lem:wpc-cut-count}.
The cut-vertex lower bound and Theorem~\ref{thm:rainbow-improve}
force equality for both rainbow parameters.
\end{proof}

\begin{corollary}
\label{cor:wpc-circular-cut-exact}
Let $G$ be a connected bagwise circular-convex WPC graph with
at least two bags and $f=0$. Suppose that, for each nonleaf
bag $B$, either
\begin{enumerate}
\item $t_B\ge2$; or
\item $P_B=\{p_B\}$ and $p_B$ belongs to every
attachment boundary of $B$.
\end{enumerate}
Then
\(
\rvc(G)=\srvc(G)=\Psi(G,T)=\ell(G).
\)
\end{corollary}

\begin{proof}
Let $\Delta_B=1$ if two attachment boundaries of $B$ are
disjoint, and $\Delta_B=0$ otherwise. The circular-convex
split-component formula gives
$\rvc(H_B)=\max\{t_B,1+\Delta_B\}$.
When $t_B\ge2$, this equals $t_B$. Under the second condition,
all attachment boundaries contain $p_B$, so $\Delta_B=0$
and $\rvc(H_B)=1=t_B$. Apply Theorem~\ref{thm:wpc-cut-exact}.
\end{proof}

% \begin{theorem}[Two-vertex hub class]\label{thm:wpc-hub-exact}
% Let $G$ be a connected WPC graph whose partition tree has
% a nonleaf bag $B_0=\{a,b\}$ of degree at least three, with
% $S_{B_0,D}=B_0$ for every neighbor $D$ of $B_0$.
% Suppose $f=0$, every other nonleaf bag $B$ satisfies
% $\rvc(H_B)=t_B$, and in each component $T_i$ of $T-B_0$
% there exists a vertex $x_i$ separated from $B_0$ by every
% cut vertex of $G$ lying in that component. Then
% \(
% \rvc(G)=\srvc(G)=\Psi(G,T)=\ell(G)+1.
% \)
% \end{theorem}

% \begin{proof}
% All markers of $H_{B_0}$ have the same neighborhood
% $\{a,b\}$; hence $\rvc(H_{B_0})=1$. Neither $a$ nor
% $b$ is a cut vertex. Lemma~\ref{lem:wpc-cut-count} and
% the local assumptions imply $\Psi(G,T)=\ell(G)+1$.

% If $\ell(G)=0$, then $G$ is not complete, so
% $\rvc(G)\ge1=\Psi(G,T)$. Otherwise suppose a rainbow
% vertex coloring uses only $\ell(G)$ colors. All cut
% vertices must have distinct colors, exhausting the palette.
% Thus the colors of $a$ and $b$ also occur at cut vertices
% $u,v$, respectively, possibly $u=v$. If $u,v$ belong to
% different components of $T-B_0$, choose the corresponding
% vertices $x_i,x_j$. If they lie in the same component, take
% $x_i$ there and $x_j$ in any other component. Every path
% between the chosen vertices passes through $u,v$ and at
% least one of $a,b$ as internal vertices. The latter vertex
% repeats the color of $u$ or $v$, contradicting the rainbow
% property. Hence $\rvc(G)\ge\ell(G)+1$, and the component-sum
% upper bound proves the claim.
% \end{proof}

\subsubsection{Case II: $f>0$ and the nonleaf subtree lies on one path}

For a leaf bag $L$ with neighbor $P(L)$, put
\(
\varepsilon_L=
\begin{cases}
1,&S_{L,P(L)}\subsetneq L,\\
0,&S_{L,P(L)}=L.
\end{cases}
\)

\begin{lemma}
\label{lem:leaf-path-distance}
Let $L,R$ be distinct leaf bags with connecting route
$L=B_0-B_1-\cdots-B_r-B_{r+1}=R$. For $1\le i\le r$,
put
\(
\delta_i=
\begin{cases}
1,&S_{B_i,B_{i-1}}\cap S_{B_i,B_{i+1}}=\varnothing,\\
0,&\text{otherwise}.
\end{cases}
\)
Choose $x\in L$ outside its attachment boundary when
$\varepsilon_L=1$, and otherwise inside it. Choose $y\in R$
analogously. Then
\(
d_G(x,y)=(r+1)+\sum_{i=1}^{r}\delta_i
                  +\varepsilon_L+\varepsilon_R.
\)
\end{lemma}

\begin{proof}
Every $x$--$y$ path uses at least $r+1$ edges between
successive bags. It additionally needs one edge inside each
intermediate bag whose incoming and outgoing boundaries are
disjoint, and one edge in each endpoint leaf when the chosen
endpoint lies outside its boundary. This gives the stated
lower bound. Since bags are cliques and the interbag edges
are complete between boundaries, all the required movements
can be realized with exactly that number of edges.
\end{proof}

\begin{theorem}\label{thm:psi-spine}
Suppose $f>0$ and there are distinct leaf bags $L,R$ with
route $L=B_0-B_1-\cdots-B_r-B_{r+1}=R$ such that
\begin{enumerate}
\item every nonleaf bag lies on this route;
\item $\rvc(H_{B_i})=1+\delta_i$ for every $1\le i\le r$;
\item $\varepsilon_L+\varepsilon_R=s+\min\{2,f\}$.
\end{enumerate}
Then
\(
\rvc(G)=\srvc(G)=\Psi(G,T)=\diam(G)-1.
\)
\end{theorem}

\begin{proof}
By the first two conditions,
\(
K(G,T)=\sum_{i=1}^{r}\rvc(H_{B_i})
      =r+\sum_{i=1}^{r}\delta_i.
\)
Using the third condition and Lemma~\ref{lem:leaf-path-distance},
\(
\Psi(G,T)=r+\sum_{i=1}^{r}\delta_i
               +\varepsilon_L+\varepsilon_R
          =d_G(x,y)-1.
\)
Hence
\(
\Psi(G,T)\le\rvc(G)\le\srvc(G)\le\Psi(G,T).
\)
All inequalities are equalities. Moreover,
$\diam(G)\ge d_G(x,y)=\Psi(G,T)+1$, whereas
$\diam(G)-1\le\rvc(G)=\Psi(G,T)$. Thus
$\diam(G)-1=\Psi(G,T)$.
\end{proof}

\begin{corollary}
\label{cor:psi-flexible}
Suppose $s=0$ and $f\ge3$. Assume that two
\emph{nontrivial} leaf bags $L,R$ with at least two attachment
vertices have a connecting route containing all nonleaf bags,
and every intermediate nonleaf bag $B_i$ satisfies
$\rvc(H_{B_i})=1+\delta_i$. Then
\(
\rvc(G)=\srvc(G)=\Psi(G,T)=K(G,T)+2.
\)
\end{corollary}

\begin{proof}
The chosen leaves have $\varepsilon_L=\varepsilon_R=1$.
Since $s=0$ and $f\ge3$, we have
$\varepsilon_L+\varepsilon_R=2=s+\min\{2,f\}$.
Apply Theorem~\ref{thm:psi-spine}.
\end{proof}

\begin{corollary}\label{cor:psi-common-spine}
Suppose the hypotheses of Corollary~\ref{cor:psi-flexible}
hold, except that its local numerical condition is replaced by:
for every nonleaf bag $B$, there is a vertex $p_B\in B$
belonging to every attachment boundary of $B$.
If $r$ is the number of nonleaf bags, then
\(
\rvc(G)=\srvc(G)=\Psi(G,T)=r+2.
\)
\end{corollary}

\begin{proof}
Every marker of $H_B$ is adjacent to $p_B$.
Thus $H_B$ has diameter two and, since $B$ is nonleaf,
$\rvc(H_B)=1$. Also the incoming and outgoing boundaries
both contain $p_B$, so $\delta_i=0$. Hence
$\rvc(H_{B_i})=1+\delta_i$ and
Corollary~\ref{cor:psi-flexible} applies, with $K(G,T)=r$.
\end{proof}

\subsubsection{Case III: $f>0$ and the nonleaf subtree branches}

In this case the original component-sum value is never attained.
More precisely, Theorem~\ref{thm:rainbow-branch-improve} gives
\(
\rvc(G)\le\srvc(G)\le\Psi(G,T)-1<\Psi(G,T).
\) In the following, we show that, under the conditions of \textbf{Case III}, this bound this also tight.

% \begin{corollary}
% \label{cor:branch-exact-certificate}
% Suppose $f>0$, $T^\circ$ branches, and
% $D(G)=\Psi(G,T)-1$. Then
% \(
% \rvc(G)=\srvc(G)=\widehat\Psi(G,T)=\Psi(G,T)-1.
% \)
% \end{corollary}

% \begin{proof}
% The lower bound $D(G)\le\rvc(G)$ and
% Theorem~\ref{thm:rainbow-branch-improve} yield
% \(
% \Psi(G,T)-1=D(G)\le\rvc(G)\le\srvc(G)
% %                           \le\Psi(G,T)-1.
% \)
% Thus equality holds throughout.
% \end{proof}

\begin{proposition}
\label{prop:branch-psi-minus-one}
For each integer $h\ge3$, there is a WPC graph $G_h$
with $f=3$, whose nonleaf subtree branches, such that
\(
\rvc(G_h)=\srvc(G_h)=\Psi(G_h,T)-1=5.
\)
\end{proposition}

\begin{proof}
First, take the central vertex $c$, three vertices $u_1,u_2,u_3$
adjacent only to $c$ outside their own branches, and three
pairwise disjoint clique bags $L_1,L_2,L_3$, each of size $h$.
In $L_i$ we choose distinct attachment vertices $a_i,b_i$;
join $u_i$ to $a_i,b_i$ and to no other vertices of $L_i$.
There are no other interbag edges. The partition tree has
routes $\{c\}-\{u_i\}-L_i$ for $i=1,2,3$. All four
nonleaf bags are singletons, so their split components
have rainbow vertex-connection number one. Hence
$K(G_h,T)=4$. Each $L_i$ is nontrivial with two attachment
vertices, so $s=0$, $f=3$, and $\Psi(G_h,T)=6$.

We choose $w_i\in L_i\setminus\{a_i,b_i\}$.
For $i\ne j$, the distance $d_G(w_i,w_j)=6$, so
$\rvc(G_h)\ge\diam(G_h)-1\ge5$.
The strict branching upper bound gives $\srvc(G_h)\le5$.
Therefore
$5\le\rvc(G_h)\le\srvc(G_h)\le5$, proving the result.
\end{proof}

\subsubsection{An improvement using recoloring}

The component-sum construction gives a strong rainbow
vertex coloring using at most $\Psi(G,T)$ colors.
However, some colors may be reused without affecting
the rainbow property of the selected shortest paths.
We now describe a polynomial-time method to reduce
the number of colors.

Let $\tau$ be a strong rainbow vertex coloring of
$G$ using a set of colors $Q$, where $q=|Q|$.
For each pair of distinct vertices $x,y$, choose
one rainbow shortest path $P_{xy}$ under $\tau$.

The \emph{interference graph} $F$ has vertex set
$Q$. Two distinct colors $a,b$ are adjacent in
$F$ if both appear on internal vertices of some
selected path $P_{xy}$.

Thus, two nonadjacent colors in $F$ can be
combined into a single color without affecting
the rainbow property of any selected shortest path.

\begin{theorem}
\label{thm:interference}
Let $F$ be the interference graph associated with
a strong rainbow vertex coloring
$\tau$ using $q$ colors. Then
\(
 \srvc(G)\le q-\nu(\overline F),
\)
where $\nu(\overline F)$ denotes the maximum
matching size in the complement of $F$.

In particular, if $\tau$ is obtained from
Theorem~\ref{thm:rainbow-improve}, then
\(
 \srvc(G)\le q-\nu(\overline F)
             \le\Psi(G,T).
\)
The recoloring can be computed in polynomial
time once the coloring and the
selected rainbow shortest paths are given.
\end{theorem}

\begin{proof}
Let $M$ be a maximum matching in $\overline F$.
For every edge $ab\in M$, we identify the colors
$a$ and $b$. Since $M$ is a matching, no color
belongs to more than one identified pair.
All unmatched colors remain distinct.
By the definition of $\overline F$, the colors
$a$ and $b$ of a matched pair never occur
together as internal colors on any selected
rainbow shortest path.
Consequently, identifying these colors cannot
create two internal vertices of the same color
on any selected path.
Since different matching edges have disjoint
endpoints, the recoloring preserves every
selected rainbow shortest path.
Each identified pair reduces the number of
colors by exactly one. The resulting coloring
therefore uses
\(
 q-|M|=q-\nu(\overline F)
\)
colors.
Hence
\(
 \srvc(G)\le q-\nu(\overline F).
\)
Finally, if the provisional coloring is obtained
from Theorem~\ref{thm:rainbow-improve}, then
$q\le\Psi(G,T)$, which gives
\(
 \srvc(G)\le q-\nu(\overline F)
             \le\Psi(G,T).
\)
\end{proof}

The recoloring bound may also certify optimality
when it coincides with a lower bound on the
rainbow vertex-connection number.

\begin{corollary}
\label{cor:matching-exact}
Let $G$ be a connected WPC graph, and let
$F$ be the interference graph associated with
a strong rainbow vertex coloring
using $q$ colors.
If
\(
 D(G)=q-\nu(\overline F),
\)
then
\(
 \rvc(G)=\srvc(G)
        =q-\nu(\overline F)=D(G).
\)
\end{corollary}

\begin{proof}
The diameter and cut-vertex lower bounds imply
\(
 D(G)\le\rvc(G).
\)
By Theorem~\ref{thm:interference},
\(
 \srvc(G)\le q-\nu(\overline F).
\)
Therefore,
\(
 D(G)\le\rvc(G)\le\srvc(G)
      \le q-\nu(\overline F).
\)
Under the stated condition, the first and
last quantities are equal. Hence all
inequalities are equalities.
\end{proof}

\section{Conclusion and Future Work}
\label{sec:conclusion}

In this paper, we studied conflict-free open- and
closed-neighborhood coloring and rainbow vertex
connection on split, starlike, and well-partitioned
chordal (WPC) graphs. For CFON and its partial
variant, we obtained polynomial-time constructive
upper bounds on circular-convex split and starlike
graphs and extended these bounds to circular-convex
WPC graphs using partition-tree structures.
For CFCN, we established an exact characterization for both
the partial and full variant on starlike graphs, proved
upper bounds of two and three colors for partial
and full CFCN on WPC graphs, respectively, and
established NP-completeness of the corresponding
one-color and two-color decision problems.
For rainbow vertex connection, we obtained exact
results on circular-convex split and starlike graphs
and got upper bounds on WPC graphs by reusing
colors across leaf bags. We also identified
subclasses of WPC graphs where these bounds
are attained.

Several questions remain open. In particular,
determining the exact complexity of CFON on
circular-convex split, starlike and WPC graphs, would be
interesting. Another direction is to develop
exact polynomial-time algorithms for rainbow
vertex connection and strong rainbow vertex
connection on broader subclasses of WPC graphs.
Finally, extending the proposed techniques to
other coloring problems and graph classes is
a natural direction for future research.

\end{document}